\documentclass[12pt]{amsart}

\usepackage{graphicx}
\usepackage{amsmath}
\usepackage{amssymb}
\usepackage{amsfonts}
\usepackage{verbatim}
\usepackage{scalefnt}
\usepackage{multirow}
\usepackage{enumerate}
\usepackage{mathtools}
\usepackage{float}
\usepackage{enumitem}
\restylefloat{figure}
\usepackage[margin=3cm]{geometry}

\usepackage{fancyhdr}
\usepackage{hyperref}

\newcommand{\F}{{\mathbb F}}
\newcommand{\Z}{{\mathbb Z}}

\newtheorem{theorem}{Theorem}

\newtheorem{proposition}[theorem]{Proposition}
\newtheorem{corollary}[theorem]{Corollary}

\newtheorem*{definition*}{Definition}

\numberwithin{equation}{section}
\numberwithin{theorem}{section}

\title[Sets avoiding six-term arithmetic progressions in $\mathbb{Z}_6^n$ are exponentially small]%
  {Sets avoiding six-term arithmetic progressions in $\mathbb{Z}_6^n$ are exponentially small}

\author{P\'eter P\'al Pach}
\email{ppp@cs.bme.hu}
\address{MTA-BME Lend\"ulet Arithmetic Combinatorics Research Group, Department of Computer Science and Information Theory, Budapest
  University of Technology and Economics, 1117 Budapest, Magyar tud\'osok
  k\"or\'utja 2., Hungary}

\author{Rich\'ard Palincza}
\email{pricsi@cs.bme.hu}
\address{MTA-BME Lend\"ulet Arithmetic Combinatorics Research Group, Department of Computer Science and Information Theory, Budapest
  University of Technology and Economics, 1117 Budapest, Magyar tud\'osok
 k\"or\'utja 2., Hungary}

\thanks{}

\begin{document}

\begin{abstract}	
We show that sets avoiding 6-term arithmetic progressions in $\mathbb{Z}_6^n$ have size at most $5.709^n$. It is also pointed out that the ``product construction'' does not work in this setting, specially, we show that for the extremal sizes in small dimensions we have $r_6(\mathbb{Z}_6)=5$, $r_6(\mathbb{Z}_6^2)=25$ and $ 116\leq r_6(\mathbb{Z}_6^3)\leq 124$.	
\end{abstract}

\date{\today}
\maketitle

\section{Introduction}

%
%

There has been great interest in finding progression-free sets in $\Z_m^n:=(\Z/(m \Z))^n$, especially when $m=3$ or $4$. Let $r_k(\mathbb{Z}_m^n)$ denote the maximal size of a set $A\subset \mathbb{Z}_m^n$
with no $k$ distinct elements in arithmetic progression.
Note that for $m=3,4,5$ the properties
``no arithmetic progression of length $3$ modulo $m$'' and ``no $3$ points on  any line'' are equivalent. The last property is also well known under the name caps.

The following is known \cite{Edel:2004,EllenbergandGijswijt:2016,Elsholtz-Pach,Croot-Lev-Pach} for the cases $k=3$, $m\in\{3,4\}$: 
\[ 2.21738\ldots ^n \leq r_3(\Z_3^n) \leq 2.755\ldots ^n, 
\]
\[ 3^n/\sqrt{n} \ll r_3(\Z_4^n) \leq 3.61\ldots^n, \]
and more generally, for primes $p\geq 3$ and some positive constant $\delta_p$
\[  r_3(\Z_p^n) \leq (p-\delta_p)^n.\]
(Note that the lower bound for $r_3(\mathbb{Z}_3^n)$ holds only for sufficiently large values of $n$, the upper bounds hold for every $n$.)
Indeed the argument yields \cite{Blasiak-Church-etal} the bound $$r_3(\Z_p^n) \leq (J(p)p)^n, $$
where 
\begin{equation}\label{jpdef}
J(p)=\frac{1}{p}\min\limits_{0<t<1} \frac{1-t^p}{(1-t)\, t^{(p-1)/3}}.
\end{equation}
As $J(p)$ is decreasing and $J(3)\leq 0.9184$ one can conclude that for every $m \geq 3$ the following holds
(see e.g. \cite{Blasiak-Church-etal} and \cite{Petrov-and-Pohoato:2018}):
\begin{equation}\label{consequence-for-all-m}
    r_3(\mathbb{Z}_m^n) \leq (0.9184m)^n 
\end{equation}
for every $m\geq 3$.

Note that the method could be applied for any finite field $\mathbb{F}_q$ with $q=p^\alpha$, however, since $r_3(\mathbb{F}_q^n)=r_3(\mathbb{F}_p^{\alpha n})$ the relevant cases are those when the prime power $q$ is a prime. (The resulting upper bound from the application to $\mathbb{F}_{p^\alpha}$ is worse than the bound coming from the case of $\mathbb{F}_p$.)


It is easy to see that the sequence $(r_3(\mathbb{Z}_m^n))^{1/n}$ converges to some limit $\alpha_{3,m}$. The main idea behind this observation is that with the help of the product construction one can bubble up constructions found in small dimensions. Namely, if $A$ avoids 3AP's in dimension $n$, then the $t$-fold direct product $\underbrace{A\times A\times \dots\times A}_t$ also avoids 3AP's in dimension $tn$. 

As $\alpha_{3,m}<m$ we may say that 3AP-free sets in $\mathbb{Z}_m^n$ are exponentially small when $m\geq 3$. Prior to this work for longer progressions it has not yet been decided in {\it any} of the cases $4\leq k\leq m$ whether $r_k(\mathbb{Z}_m^n)$ is also exponentially small or of order of magnitude $(m-o(1))^n$ (as $n\to \infty$). In this note we will prove that whenever $6\mid m$ and $k\in \{4,5,6\}$ the quantity $r_k(\mathbb{Z}_m^n)$ is exponentially small, specially, $r_6(\mathbb{Z}_6^n)\leq 5.709^n$. It is tempting to also formulate this statement as $\lim (r_6(\mathbb{Z}_6^n))^{1/n}\leq 5.709$, however, somewhat surprisingly, we do not see a proof of the statement that 
$r_6(\mathbb{Z}_6^n)^{1/n}$ converges (although we believe it surely does). The convergence is not immediate, because the product construction does not work in general. When $k=3$ or $m$ is a prime power, the $t$-fold direct product $\underbrace{A\times A\times \dots \times A}_t$ avoids $k$-AP's when $A$ itself is $k$-AP-free, however, for general $k$ and $m$ this fails to hold.  
Let us illustrate this by the case $k=6,m=6$. In dimension 1 we clearly have $r_6(\mathbb{Z}_6)=5$, and, for instance, the set $A=\{0,1,2,3,4\}$ is 6AP-free. By taking $A\times A=\{0,1,2,3,4\}\times \{0,1,2,3,4\}$ we obtain a 25-element subset of $\mathbb{Z}_6^2$ which contains the following 6AP:
$$(0,0),(2,3),(4,0),(0,3),(2,0),(4,3).$$
Although the product construction is not applicable, the value of $r_6(\mathbb{Z}_6^2)$ still turns out to be $25=5^2$, however, we will show that $r_6(\mathbb{Z}_6^3)<125=(r_6(\mathbb{Z}_6))^3$.

Summarizing our results we prove the following bounds:

\begin{theorem}\label{6APthm}
For sets without arithmetic progression of length $6$
we have the following results in small dimensions:
\[r_6(\Z_6^1)=5, r_6(\Z_6^2)=25, 116\leq r_6(\Z_6^3)\leq 124.\]
\end{theorem}

\begin{theorem}\label{thm-6AP-exp-small}
For sets without arithmetic progression of length $6$
we have the following results:
$$4.434^n\leq 2^nr_3(\mathbb{Z}_3^n)\leq r_6(\mathbb{Z}_6^n)\leq 5.709^n,$$
assuming that $n$ is sufficiently large.
\end{theorem}

If $6\mid m$, then $\mathbb{Z}_6^n$ is a subgroup of $\mathbb{Z}_m^n$, and by using the bound from Theorem~\ref{thm-6AP-exp-small} in each of the $(m/6)^n$ cosets the following corollary is obtained: 

\begin{corollary}
If $6\mid m$ and $k\in \{4,5,6\}$, then $r_k(\mathbb{Z}_m^n)\leq (0.948m)^n,$ if $n$ is sufficiently large.
\end{corollary}

Finally, we provide another upper bound for $r_6(\mathbb{Z}_6^n)$ in terms of $r_3(\mathbb{Z}_3^n)$. 

\begin{theorem}
\label{thm-6AP-exp-small-2}
For sets without arithmetic progression of length $6$
we have the following result:
$$r_6(\mathbb{Z}_6^n)\leq 2^{n+1}\sqrt{3^nr_3(\mathbb{Z}_3^n)}.$$
\end{theorem}
Note that by using the bound $r_3(\mathbb{Z}_3^n)\leq 2.756^n$ Theorem~\ref{thm-6AP-exp-small-2} implies that $r_6(\mathbb{Z}_6^n)\leq 5.75^n$ which bound is worse than the one in Theorem~\ref{thm-6AP-exp-small}, however, if $r_3(\mathbb{Z}_3^n)\leq 2.69^n$, then Theorem~\ref{thm-6AP-exp-small-2} gives a better estimation than Theorem~\ref{thm-6AP-exp-small}.

The paper is organized as follows:
In Section~\ref{sec-subset}  we give a reformulation for the problem of finding $r_k(\mathbb{Z}_6^n)$ with $k\in \{3,4,5,6\}$ in terms of possible total sizes of systems of subsets of $\mathbb{Z}_3^n$ satisfying certain properties. In Section~\ref{sec-proofs} we prove Theorem~\ref{6APthm}, Theorem~\ref{thm-6AP-exp-small} and Theorem~\ref{thm-6AP-exp-small-2}. Some concluding remarks are given in Section~\ref{sec-concl}.

\section{Subset reformulation}\label{sec-reformulation}\label{sec-subset}

We may express $\mathbb{Z}_6^n$ as $\mathbb{Z}_6^n=F\oplus R$, where $F=\{0,2,4\}\cong \mathbb{Z}_3^n$ and $R=\{0,3\}^n\cong \mathbb{Z}_2^n$. A sequence $a_1=f_1+r_1,a_2=f_2+r_2,\dots,a_k=f_k+r_k$ (where $f_i\in F, r_i\in R$) forms an arithmetic progression in $\mathbb{Z}_6^n$ if and only if $f_1,f_2,\dots,f_k$ is an arithmetic progression in $\mathbb{Z}_3^n$ and $r_1,r_2,\dots,r_k$ is an arithmetic progression in $\mathbb{Z}_2^n$, respectively. Note that if the elements are distinct, then $k\leq 6$. If $k=3$, then the progression consists of pairwise different elements if and only if $f_1,f_2,f_3$ are distinct. Since the sequence $r_1,r_2,\dots$ is alternating, for $k\in \{4,5,6\}$ the necessary and sufficient conditions for getting $k$ distinct elements is that $f_1,f_2,f_3$ are distinct and $r_1,r_2$ are distinct. Using this decomposition we may reformulate the property that ``a subset $A\subseteq \mathbb{Z}_6^n$ avoids $k$-term arithmetic progressions'' in terms of a property of systems of subsets of $\mathbb{Z}_3^n$. Namely, let $A(r)=\{f\in \mathbb{Z}_3^n: f+r\in A\}$ for $r\in R$ and
let us define properties $(*)_3,(*)_4,(*)_5,(*)_6$ as follows:

The system of subsets $A(r)$ ($r\in \mathbb{Z}_2^n$) satisfies 
\begin{itemize}
    \item property $(*)_3$, if $A(r')\cup A(r'')$ is 3AP-free for every pair $r',r''\in \mathbb{Z}_2^n$,
    \item property $(*)_4$, if it is not possible to choose two different indices $r',r''\in \mathbb{Z}_2^n$ and a 3AP $a,b,c$ in $\mathbb{Z}_3^n$ such that $a,b\in A(r')$ and $a,c\in A(r'')$,
    \item property $(*)_5$, if it is not possible to choose two different indices $r',r''\in \mathbb{Z}_2^n$ and a 3AP $a,b,c$ in $\mathbb{Z}_3^n$ such that $a,b,c\in A(r')$ and $a,b\in A(r'')$,
    \item property $(*)_6$, if $A(r')\cap A(r'')$ is 3AP-free for every pair of distinct indices $r',r''\in \mathbb{Z}_2^n$.
\end{itemize}
Note that in this reformulation $\mathbb{Z}_2^n$ serves only as an index set of size $2^n$, its structure does not play any role.

Let us summarize in  the following statement how the reformulation can be used to study the $r_k(\mathbb{Z}_6^n)$ values.

\begin{proposition}
Let $k\in \{3,4,5,6\}$. The maximum total size of a system of subsets $A(r)\subseteq \mathbb{Z}_3^n$ ($r\in \mathbb{Z}_2^n$) satisfying property $(*)_k$ is $r_k(\mathbb{Z}_6^n)$.
\end{proposition}
\begin{proof}
The statements immediately follow from the structural description of arithmetic progressions in $\mathbb{Z}_6^n$.
\end{proof}

Let us mention that the problem of determining the size of the largest 3AP-free subset of $\mathbb{Z}_6^n$ is equivalent with doing so in case of $\mathbb{Z}_3^n$:

\begin{proposition}\label{prop-3AP-mod6} For sets without arithmetic progression of length three  the following holds:
$$r_3(\mathbb{Z}_6^n)=2^nr_3(\mathbb{Z}_3^n).$$
\end{proposition}
\begin{proof}
If $A_0\subseteq \mathbb{Z}_3^n$ is 3AP-free, then the system $A(x)=A_0$ ($x\in\mathbb{Z}_2^n$) satisfies property $(*)_3$, thus $r_3(\mathbb{Z}_6^n)\geq 2^nr_3(\mathbb{Z}_3^n)$.

On the other hand, if $\sum\limits_{x\in\mathbb{Z}_2^n} |A(x)|>2^nr_3(\mathbb{Z}_3^n)$, then for some $x$ we have $|A(x)|>r_3(\mathbb{Z}_3^n)$, thus $A(x)$ contains a 3AP, and $(*)_3$ fails to hold. Hence, $r_3(\mathbb{Z}_6^n)=2^nr_3(\mathbb{Z}_3^n)$.
\end{proof}

In fact the argument only used that 6 has residue 2 modulo 4, and in general it yields the following statement:
\begin{proposition} If $m=4M+2$ for some integer $M$, then
$$r_3(\mathbb{Z}_m^n)=2^nr_3(\mathbb{Z}_{m/2}^n).$$
\end{proposition}
While studying $r_3(\mathbb{Z}_m^n)$ there are some technical differences between the cases when $m$ is odd and when $m$ is divisible by 4, but the case when $m$ is an even number not divisible by 4 simply reduces to the odd case. We shall mention that for certain composite values of $m$ there has been some improvements on the trivial corollaries of the prime case, like $r_3(\mathbb{Z}_9^n)\leq 3^nr_3(\mathbb{Z}_3^n)$. Namely, the method was adapted to odd prime powers \cite{Blasiak-Church-etal,Petrov,Speyer} and also to the technically more difficult even case for $m=2^3=8$. \cite{Petrov-and-Pohoato:2018}

\section{Proofs} \label{sec-proofs}

\begin{proof}[Proof of Theorem~\ref{6APthm}]

{\bf Dimension 1.}
Clearly, $r_6(\mathbb{Z}_6^1)=5$. Any 5-element subset of $\mathbb{Z}_6^1$ is trivially 6AP-free.

{\bf Dimension 2.}
Now, we show that $r_6(\mathbb{Z}_6^2)=25$. Using the reformulation from Section~\ref{sec-reformulation} we are interested in the maximal possible total size of a system of four subsets of $\mathbb{Z}_3^2$ satisfying property $(*)_6$. That is, we would like to determine the maximum of
$\sum\limits_{i=1}^4 |A_i|$, where $A_i\subseteq \mathbb{Z}_3^2\ (1\leq i\leq 4)$ such that no 3AP is contained in at least two of the subsets $A_i$. The total number of 3AP's in $\mathbb{Z}_3^2$ is $\frac{9\cdot8}{6}=12$, thus the four subsets $A_1,\ A_2,\ A_3,\ A_4$ can contain at most twelve 3AP's in total. It is easy to determine the smallest possible number of 3AP's that must be contained in a subset of a given size (by hand or by a computer search). Let us summarize the results in the table below:
\medskip
\begin{center}
\begin{tabular}{ |c||c|c|c|c|c|c|c|c|c|c| } 
\hline
size of $A$ & 0 & 1 & 2 & 3 & 4 & 5 & 6 & 7 & 8 & 9\\
\hline
min \#3AP in $A$ & 0 & 0 & 0 & 0 & 0 & 1 & 2 & 5 & 8 & 12 \\
\hline
\end{tabular}
\end{center}
\medskip

Let $x_i$ denote the number of $i$-element subsets among $A_1,A_2,A_3,A_4$ (where $0\leq i \leq 9$).
Since each 3AP can appear in at most one set $A_i$, the optimal value for $\sum\limits_{i=1}^4 |A_i|$ can not be more than the solution of the following integer program:

\medskip
\centerline{$\begin{array}{l}
\max x_1+2x_2+3x_3+4x_4+5x_5+6x_6+7x_7+8x_8+9x_9\\
\hbox{subject to}\\
x_0+x_1+x_2+x_3+x_4+x_5+x_6+x_7+x_8+x_9=4\\
x_5+2x_6+5x_7+8x_8+12x_9\le 12\\
x_0,x_1,x_2,x_3,x_4,x_5,x_6,x_7,x_8,x_9: \text{nonnegative integers}\\
\end{array}$}
\medskip

(The first constraint ensures that four subsets are chosen, and the second constraint holds, since the total number of 3AP's contained in the four subsets can not be more than the total number of 3AP's in $\mathbb{Z}_3^2$.)

By solving the above integer program we obtain that the optimal value is 25 which is attained at $x_6=3,x_7=1$ (everything else is 0). That is, to achieve 25, one of the subsets must have size 7, and the three other subsets must have size 6.

By symmetry, we may assume that $A_1=\mathbb{Z}_3^2\setminus \{u,v\}$, where $u$ and $v$ are two different elements. Let $w=-u-v$ be the third point on the line $uv$. Let $\alpha$ denote the direction of the line $uv$. Note that in $\mathbb{Z}_3^2$ there are four possible directions, let us denote the other three directions by $\beta,\gamma$ and $\delta$.

Note that $A_2,A_3,A_4$ must have size 6 and each of them must contain exactly two 3AP's. In $\mathbb{Z}_3^2$ there are two types of 6-element sets: the complement of a 6-element set is either an affine line or not. To contain only two 3AP's the sets $A_2,A_3,A_4$ must all be the complements of affine lines, in other words, each of them is a union of two parallel lines. Moreover, these lines must not be parallel with the line $uv$, otherwise at least one of them would be contained in two subsets (in $A_1$ and here).

Also, none of these lines can go through  $w$, as this would result in a 3AP contained both in $A_1$ and here.

Finally, a line from $A_i$ and a line from $A_j$ (where $2\leq i<j\leq 4$) must not be parallel with each other because of similar reasons. That is, we may assume that $A_2,A_3,A_4$ are the unions of two-two lines of directions $\beta, \gamma, \delta$, respectively.

Therefore, $A_2,A_3,A_4$ can be characterized as follows: $A_2,A_3,A_4$ are all the unions of two parallel lines, where the directions of the lines are $\beta,\gamma,\delta$ resp., furthermore each line goes through $u$ or $v$. 
(Thus $\{A_2,A_3,A_4\}$ is uniquely determined.)

The obtained system $\{A_1,A_2,A_3,A_4\}$ satisfies the conditions, since:
\begin{itemize}
    \item $A_1$ contains two 3AP's with direction $\alpha$ and three more 3AP's that contain $w$.
    \item None of the 3AP's contained in $A_2,A_3,A_4$ have direction $\alpha$ and none of them contains $w$.
    \item The two-two lines contained in $A_2,A_3,A_4$ have directions $\beta,\gamma,\delta$, respectively.
\end{itemize}

Hence, we proved that the largest 6AP-free set in $\mathbb{Z}_6^2$ has size 25 (and it is unique in the above described sense).

\medskip

{\bf Dimension 3.}
Analogously to the previous case, with a quick computer check we find that the minimum number of 3AP's 
that must be contained in subsets of $\mathbb{Z}_3^3$ of given sizes are the numbers below. (Let $m_j$ denote the minimum number of 3AP's that must be contained in a set of size $j$.)
\medskip
\begin{center}
\begin{tabular}{ |c||c|c|c|c|c|c|c|c|c|c|c|c|c|c| } 
\hline
size of $A$  ($j$) & 0 & 1 & 2 & 3 & 4 & 5 & 6 & 7 & 8 & 9 & 10 & 11 & 12 & 13\\
\hline
min \#3AP in $A$ ($m_j$) & 0 & 0 & 0 & 0 & 0 & 0 & 0 & 0 & 0 & 0 & 2 & 3 & 4 & 7 \\
\hline
\end{tabular}
\end{center}
\medskip
\begin{center}
\begin{tabular}{ |c||c|c|c|c|c|c|c|c|c|c|c|c|c|c| } 
\hline
size of $A$ ($j$) & 14 & 15 & 16 & 17 & 18 & 19 & 20 & 21 & 22 & 23 & 24 & 25 & 26 & 27\\
\hline
min \#3AP in $A$ ($m_j$) & 10 & 13 & 16 & 20 & 24 & 33 & 42 & 51 & 60 & 70 & 80 & 92 & 104 & 117 \\
\hline
\end{tabular}
\end{center}
\medskip

Let $x_i$ denote the number of $i$-element subsets among $A_1-A_8$ (where $0\leq i \leq 27$). 

Since each 3AP can appear in at most one set $A_i$, the optimal value for $\sum\limits_{i=1}^8 |A_i|$ can not be more than the solution of the following integer program:

\medskip
\centerline{$\begin{array}{l}
\max \sum\limits_{i=1}^{27} ix_i\\
\hbox{subject to}\\
\sum\limits_{i=0}^{27} x_i=8\\
\sum\limits_{i=0}^{27}m_ix_i\le 117\\
x_0,x_1,\dots,x_{27}: \text{nonnegative integers}\\
\end{array}$}
\medskip

With the help of an IP solver we obtained that the optimum is 124 yielding the bound
$$r_6(\mathbb{Z}_6^3)\leq 124.$$

Turning to the lower bound, with computer help we found the following  construction where the total size of the eight subsets is 116:

\bigskip

\hspace{-0.6cm}\begin{minipage}[t]{165mm}
\begin{tabular}{|c|c|c|}
\hline
&&\\ \hline
o&o&o\\ \hline
o&o&o\\ \hline \hline
o&&\\ \hline 
&o&\\ \hline
o&o&o\\ \hline \hline
&&\\ \hline
o&o&\\ \hline
&o&o\\ \hline
\multicolumn{3}{c}{\rule{0pt}{3ex}$A_1$}

\end{tabular}\ \ \ 
\begin{tabular}{|c|c|c|}
\hline
o&o&\\ \hline
o&&o\\ \hline
&o&o\\ \hline \hline
&o&o\\ \hline 
&o&\\ \hline
&&o\\ \hline \hline
o&&\\ \hline
&o&o\\ \hline
o&o&\\ \hline
\multicolumn{3}{c}{\rule{0pt}{3ex}$A_2$}

\end{tabular} \ \ \ 
\begin{tabular}{|c|c|c|}
\hline
o&&o\\ \hline
&o&\\ \hline
o&o&\\ \hline \hline
&&o\\ \hline 
o&&o\\ \hline
o&o&\\ \hline \hline
o&&o\\ \hline
&o&o\\ \hline
o&&\\ \hline
\multicolumn{3}{c}{\rule{0pt}{3ex}$A_3$}

\end{tabular} \ \ \ 
\begin{tabular}{|c|c|c|}
\hline
&o&o\\ \hline
o&o&\\ \hline
&&\\ \hline \hline
o&o&\\ \hline 
o&&o\\ \hline
&o&o\\ \hline \hline
o&&o\\ \hline
o&&\\ \hline
o&o&\\ \hline
\multicolumn{3}{c}{\rule{0pt}{3ex}$A_4$}

\end{tabular} \ \ \ 
\begin{tabular}{|c|c|c|}
\hline
o&&\\ \hline
&&\\ \hline
&&\\ \hline \hline
o&o&o\\ \hline 
o&o&o\\ \hline
o&&\\ \hline \hline
o&o&o\\ \hline
o&o&o\\ \hline
&&\\ \hline
\multicolumn{3}{c}{\rule{0pt}{3ex}$A_5$}

\end{tabular} \ \ \ 
\begin{tabular}{|c|c|c|}
\hline
\phantom{o}&o&\\ \hline
&o&o\\ \hline
&o&\\ \hline \hline
&o&o\\ \hline 
&o&o\\ \hline
&o&o\\ \hline \hline
&o&o\\ \hline
&&o\\ \hline
&&o\\ \hline
\multicolumn{3}{c}{\rule{0pt}{3ex}$A_6$}

\end{tabular} \ \ \ 
\begin{tabular}{|c|c|c|}
\hline
&o&o\\ \hline
o&&o\\ \hline
&o&\\ \hline \hline
o&&\\ \hline 
o&&\\ \hline
&o&o\\ \hline \hline
o&o&\\ \hline
&o&o\\ \hline
&&o\\ \hline
\multicolumn{3}{c}{\rule{0pt}{3ex}$A_7$}

\end{tabular} \ \ \ 
\begin{tabular}{|c|c|c|}
\hline
o&&o\\ \hline
o&&\\ \hline
o&&o\\ \hline \hline
&o&\\ \hline 
&o&o\\ \hline
&&\\ \hline \hline
&o&\\ \hline
o&o&\\ \hline
o&o&o\\ \hline
\multicolumn{3}{c}{\rule{0pt}{3ex}$A_8$}

\end{tabular}
\end{minipage}


\bigskip

Hence,
$$116\leq r_6(\mathbb{Z}_6^3).$$

\end{proof}


\begin{proof}[Proof of Theorem~\ref{thm-6AP-exp-small}]

The lower bound follows from Proposition~\ref{prop-3AP-mod6} and Edel's \cite{Edel:2004} lower bound for $r_3(\mathbb{Z}_3^n)$.


For proving the upper bound 
it suffices to show that $\sum \limits_{i\in I}|A_i|\leq 5.709^n$ if the system of subsets (of $\mathbb{Z}_3^n$) $\{A_i:i\in I\}$ satisfies property $(*)_6$ and $|I|=2^n$.

We will use a supersaturation extension \cite{FL} of the cap set result. This says that any subset of $\mathbb{Z}_3^n$ of density $\alpha$
 has three-term arithmetic progression density at least $\alpha^C$, where $C \approx 13.901$ is an
explicit constant \footnote{Namely, $C=1+\frac{\log 3}{\log(3/\alpha)}$, where $\alpha=3J(3)=2.755\dots $}. (Note that this includes counting trivial three-term arithmetic progressions.)

Let $\beta=3/2^{1/C}\approx 2.854$, then we have $\beta^C=\frac{3^C}{2}$. 
The total size of subsets having size at most $\beta^n$ is at most $2^n\beta^n$. Now, we consider the subsets with size larger than $\beta^n$. Let $m_i$ denote the number of those subsets whose size lies in $(2^i\beta^n,2^{i+1}\beta^n]$. Since each 3AP can occur in at most one set, we obtain that
$$m_i (2^i\beta^n/3^n)^C\leq 1$$
yielding that $m_i\leq (3/\beta)^{Cn}2^{-iC}$. Therefore, the total size of subsets of size larger than $\beta^n$ is at most $$\sum\limits_{i=0}^{\infty} m_i2^{i+1}\beta^n\leq \sum\limits_{i=0}^{\infty} (3/\beta)^{Cn}2^{-iC}2^{i+1}\beta^n =(2\beta)^n \sum\limits_{i=0}^{\infty} 2^{1-(C-1)i}\leq 2.001 (2\beta)^n.  $$
Hence, by adding up the obtained upper bounds for sets of size at most $\beta^n$ and larger than $\beta^n$ it is obtained that $\sum |A_i|\leq 3.001(2\beta)^n$. 

\end{proof}



\begin{proof}[Proof of Theorem~\ref{thm-6AP-exp-small-2}]
It suffices to prove that $S:=\sum \limits_{i\in I}|A_i|\leq 2^{n+1}\sqrt{3^nr_3(\mathbb{Z}_3^n)}$ if the system of subsets (of $\mathbb{Z}_3^n$) $\{A_i:i\in I\}$ satisfies property $(*)_6$ and $|I|=2^n$.

Let us enumerate the elements of $\mathbb{Z}_3^n$ by the positive integers from $[3^n]$. For $i\in I$ let $v_i$ be the characteristic vector of $A_i$, that is, the $j$th entry of $v_i$ is 1 if the element (from $\mathbb{Z}_3^n$) labeled by $j$ is contained in $A_i$ and 0 otherwise. Let $w:=\sum\limits_{i\in I} v_i$, denote the entries of $w$ by $w_1,\dots,w_{3^n}$. Note that $w_1+\dots+w_{3^n}=\sum\limits_{i\in I} |A_i|=S$.

By the Cauchy inequality
\begin{equation}\label{eq_1}
w^2=w_1^2+\dots+w_{3^n}^2\geq \frac{(w_1+\dots+w_{3^n})^2}{3^n}=\frac{S^2}{3^n}.
\end{equation}

Since $A_i\cap A_j$ is 3AP-free for any two different indices $i,j\in I$ we have $v_iv_j\leq r_3(\mathbb{Z}_3^n)$. Therefore,
\begin{equation}\label{eq_2}
w^2=\sum\limits_{i\in I} v_i^2+\sum\limits_{i,j\in I, i\ne j} v_iv_j\leq S+2^{2n}r_3(\mathbb{Z}_3^n).
\end{equation}
By comparing \eqref{eq_1} and \eqref{eq_2} we obtain that $S^2-3^nS-2^{2n}3^nr_3(\mathbb{Z}_3^n)\leq 0$ which yields 
$$S\leq \frac{3^n+\sqrt{3^{2n}+2^{2n+2}3^nr_3(\mathbb{Z}_3^n)}}{2}<2^{n+1}\sqrt{3^nr_3(\mathbb{Z}_3^n)}.$$

\end{proof}

\section{Concluding remarks} \label{sec-concl}

In this paper we prove that $r_6(\mathbb{Z}_6^n)\leq 5.709^n$, which implies that $r_k(\mathbb{Z}_m^n)$ is exponentially smaller than $m^n$ when $6\mid m$ and $k\in \{4,5,6\}$. Previously this was known only for the cases $3=k\leq m$, and according to our knowledge there is no pair of $k,m$ with $3\leq k\leq m$ such that $r_k(\mathbb{Z}_m^n)=(m-o(1))^n$ is known to be true. 

\section{Acknowledgements}

Both authors were supported by the Lend\"ulet program of the Hungarian Academy of Sciences (MTA). PPP was also supported by the National Research, Development and Innovation Office NKFIH (Grant Nr. K124171, K129335 and BME NC TKP2020). RP was also supported by  the BME-Artificial Intelligence FIKP grant of EMMI (BME FIKP-MI/SC).

\medskip


\begin{thebibliography}{99}




\bibitem{Blasiak-Church-etal}
J.~Blasiak, T.~Church, H.~Cohn, J.~Grochow, E.~Naslund, W.~Sawin, C.~Umans, 
On cap sets and the group-theoretic approach to matrix multiplication. 
Discrete Anal. 2017, Paper No. 3, 27 pp. 

\bibitem{Croot-Lev-Pach}
E.~Croot, V.F.~Lev, P.P.~Pach,
Progression-free sets in $\Z_4^n$ are exponentially small, 
Ann. of Math. (2) 185 (2017), no. 1, 331--337. 

\bibitem{Edel:2004}
Y. Edel, {Extensions of generalized product caps}, Des. Codes
Cryptography
  \textbf{31} (2004), 5 -- 14.

\bibitem{EllenbergandGijswijt:2016}
J.~S.~Ellenberg, D.~Gijswijt, 
On large subsets of $\F_q^n$ with no three-ter arithmetic progression. Ann. of Math. (2) 185 (2017), no. 1, 339--343. 


\bibitem{Elsholtz-Pach} C.~Elsholtz, P.~P.~Pach, Des. Codes
Cryptography (2020),  https://doi.org/10.1007/s10623-020-00769-0



\bibitem{FL} J.~Fox and L.~M.~Lov\'asz,
\newblock{A tight bound for Green’s arithmetic triangle removal lemma in vector
spaces},
\newblock Advances {\bf 321} (2017)  287--297.

\bibitem{Petrov}  F.~Petrov, Combinatorial results implied by many zero divisors in a group ring, arXiv: 1606.03256

\bibitem{Petrov-and-Pohoato:2018} 
F.~Petrov, C.~Pohoata,
 Improved Bounds for Progression-Free Sets in $C_8^n$, 
Israel J. Math. 236 (2020), no. 1, 345--363.

\bibitem{Speyer} D.~Speyer, https://sbseminar.wordpress.com/2016/07/08/bounds-for-sum-free-sets-inprime-power-cyclic-groups-three-ways/
\end{thebibliography}
\end{document}